\documentclass{amsart}
\usepackage{amsmath,amssymb,enumitem}
\usepackage[utf8]{inputenc}
\usepackage{amsmath, amssymb, amsfonts}
\usepackage{geometry}
\usepackage{amsthm}

\usepackage{xcolor}
\usepackage{setspace}
\usepackage{graphicx}
\usepackage{caption}

\setlist[enumerate]{label=(\arabic*)}

\newtheorem{theorem}{Theorem}[section]
\newtheorem{lemma}[theorem]{Lemma}
\newtheorem{definition}[theorem]{Definition}

\theoremstyle{definition}

\usepackage{amsthm}

\newtheorem*{question}{Question}

\theoremstyle{plain}
\newtheorem{remark}[theorem]{Remark}

\numberwithin{equation}{section}

\title[A positive answer to Barsotti and Carman's question on unit groups of Burnside rings. ]{A positive answer to Barsotti and Carman's question on unit groups of Burnside rings of fusion systems}

\makeatletter
\def\@settitle{
  \begin{center}
  \baselineskip14\p@\relax
  \bfseries
  \@title
  \end{center}
}
\makeatother

\author{Chunbo Zhou}
\address{Chunbo Zhou, School of Mathematics and Statistics, Hubei University, Wuhan, Hubei, 430062, P.~R.~China}
\email{chunbozhou6@163.com}

\author{Heguo Liu}
\address{Heguo Liu, School of Mathematics and Statistics, Hainan University, Haikou, Hainan, 570228, P.~R.~China}
\email{ghliu@hainanu.edu.cn}

\renewenvironment{abstract}
{
    \par\noindent\textbf{Abstract.}\quad\ignorespaces
}
{
    \par\ignorespacesafterend\vspace{1em}
}

\makeatletter
\def\@setauthors{
  \begingroup
  \def\thanks{\protect\thanks@warning}
  \trivlist
  \centering\footnotesize \@topsep30\p@\relax
  \advance\@topsep by -\baselineskip
  \item\relax
  \author@andify\authors
  \def\\{\protect\linebreak}
  \authors
  \ifx\@empty\contribs
  \else
    ,\penalty-3 \space \@setcontribs
    \@closetoccontribs
  \fi
  \endtrivlist
  \endgroup
}

\def\section{\@startsection{section}{1}
  \z@{.7\linespacing\@plus\linespacing}{.5\linespacing}
  {\normalfont\bfseries\centering}}
\def\subsection{\@startsection{subsection}{2}
  \z@{.5\linespacing\@plus.7\linespacing}{.5\linespacing}
  {\normalfont\bfseries}}
\makeatother

\begin{document}
\maketitle
\begin{abstract}
Let $\mathcal F$ be a saturated fusion system on a finite $2$-group $S$, and let $B(\mathcal F)$ be the Burnside ring of $\mathcal F$. Barsotti and Carman(\cite{A1})proved that, when $S$ is abelian,
$$ B(\mathcal F)^\times=B(S)^\times \quad\Longleftrightarrow\quad \mathcal F=\mathcal F_S(S). $$
Barsotti and Carman asked on page $728$ of \cite{A1} whether there exists a counterexample to Corollary 9.3 for a nonabelian $2$-group. We have given a positive answer. Our main conclusion is that the following are equivalent for all finite $2$-groups:
\begin{enumerate}
    \item $ \mathcal F = \mathcal F_S(S)$;
    \item $ B(\mathcal F)^\times = B(S)^\times $;
    \item For every maximal subgroup $M$ of $S$, $1-[S/M] \in B(\mathcal F)^\times$.
\end{enumerate}
  \medskip

\end{abstract}


\section{Introduction}
For a saturated fusion system on a finite $2$-group $S$, we study the unit group of the Burnside ring of the saturated fusion system $B(\mathcal F)^\times$, as defined by Matthew Gelvin and Sune Reeh. Indeed, $B(S)^\times$ is an elementary abelian $2$-group, and hence $B(\mathcal F)^{\times} \leq B(S)^{\times}$ is naturally an elementary abelian $2$-group. When $S$ is abelian $2$-group, Barsotti and Carman completely determine the unit group $B(\mathcal F)^\times$ in \cite{A1}, as the following.

\begin{theorem}\cite[Corollary 9.3]{A1}\label{Corollary 9.3}
Let $S$ be an abelian $2$-group and $\mathcal F$ a saturated fusion system on $S$. Then $B(\mathcal F)^\times=B(S)^\times $ if and only if $\mathcal F=\mathcal F_S(S) $.
\end{theorem}

A natural question is whether Theorem \ref{Corollary 9.3} holds for nonabelian $2$-groups. At the end of paper \cite{A1}, Barsotti and Carman posed the following question:

\begin{question}\cite[p.728]{A1}
Does there exist a counterexample to Theorem \ref{Corollary 9.3}(\cite[Corollary 9.3]{A1}) for a nonabelian $2$-group?
\end{question}

This paper gives a positive answer to this question. Recall that for a saturated fusion system $ \mathcal F $ on a finite $p$-group $S$, A. Díaz, A. Glesser, S. Park and R. Stancu proved the following: if the restriction map $H^{1}(\mathcal F; \mathbb{F}_{p} )\rightarrow H^{1}(S; \mathbb{F}_{p})$ is an isomorphism, then $ \mathcal F = \mathcal F_S(S)$ (\cite[Corollary 1.2]{DGPS}). This is the key to the proof of the main result of this paper, as following.

\begin{theorem}\label{theorem:A}
Let $\mathcal F$ be a saturated fusion system on a finite $2$-group $S$, and let $B(\mathcal F)$ be the Burnside ring of $\mathcal F$. The following are equivalent:
\begin{enumerate}
    \item $ \mathcal F = \mathcal F_S(S)$;
    \item $ B(\mathcal F)^\times = B(S)^\times $;
    \item For every maximal subgroup $M$ of $S$, $1-[S/M] \in B(\mathcal F)^\times$.
\end{enumerate}
\end{theorem}

This paper is organized as follows. In Section $2$, we present some background material on Burnside rings and fusion systems, together with several lemmas that will be used in the proofs of the main results of this paper. In Subsection $2.3$, we introduced the most important tool for proving the main conclusion of this paper, namely Lemma \ref{thm:DGPS-cor.1.2} (\cite[Corollary 1.2]{DGPS}). In Section $3$, we give the proof of the main result of this paper.

\section{Burnside rings and fusion systems}\label{sec:prelim}

\subsection{The Burnside ring of a finite group}
Let $G$ be a finite group. Let $B_{+}(G)$ denote the set of isomorphism classes of finite $G$-sets, equipped with the additive law $+$ induced by the disjoint union of $G$-sets, and with the multiplicative law $\times$ induced by the cartesian product of $G$-sets $($with diagonal action$)$. This forms a semiring.

The Burnside ring $B(G)$ of G is defined to be the Grothendieck group of $B_{+}(G)$. It is commutative and its identity element is the class $[G/G]$ of the one-element $G$-set.
One can characterize a $G$-set up to isomorphism using the following fundamental theorem of Burnside (\cite[Chapter XII, Theorem I]{b2}):

\begin{lemma}[Burnside]\cite[Chapter XII, Theorem I]{b2}\label{lemma:Burnside}
Let $G$ be a finite group, and $X$ and $Y$ be finite $G$-sets. Then the following are equivalent:
\begin{enumerate}
    \item The $G$-sets $X$ and $Y$ are isomorphic.
    \item For any subgroup $H$ of $G$, the sets $X^H$ and $Y^H$ have the same cardinality.
\end{enumerate}
\end{lemma}

\begin{proof}
It is clear that (1) implies (2), since any $G$-set isomorphism $X \to Y$ induces a bijection $X^H \to Y^H$ on the sets of fixed points by any subgroup $H$ of $G$.

To show the converse, observe that any finite $G$-set $X$ can be written up to isomorphism as
\[
X = \bigsqcup_{K \in [Cl(G)]} a_K(X)\, G/K
\]
for some $a_K(X) \in \mathbb{N}$, where $a_K(X)\, G/K$ denotes the disjoint union of $a_K(X)$ copies of $G/K$ and $[Cl(G)]$ denotes the set of representatives of $Cl(G)$.

Now if (2) holds, for any $H \in Cl(G)$, there is an equation
\[
\sum_{K \in [Cl(G)]} \big(a_K(X) - a_K(Y)\big)\big|(G/K)^H\big| = 0
\]

The matrix $m$ of this system of equations is given by
\[
m(H,K) = \big|(G/K)^H\big| = \big|\{x \in G/K \mid H^x \subseteq K\}\big|
\]
for $K,H \in [Cl(G)]$. In particular the entry $m(H,K)$ is non-zero if and only if some conjugate of $H$ is contained in $K$.

If the set $Cl(G)$ is given a total ordering $\preceq$ such that $H \preceq K$ implies $|H| \leq |K|$, then the matrix $m$ is upper triangular, with non-zero diagonal coefficient $m(H,H) = |N_G(H) : H|$. In particular $m$ is non-singular, and it follows that $a_K(X) = a_K(Y)$, for any $K \in Cl(G)$, and the $G$-sets $X$ and $Y$ are isomorphic.
\end{proof}

\begin{lemma}\cite[Formula(3.1)]{A2}\cite[Lemma 1.1]{A9}
Let $P$ and $Q$ be subgroups of $G$, and let $P\backslash G/Q$ denote the set of double cosets $\overline{g}:=PgQ$ for $g\in G$. We have:
\[
[G/P]\times [G/Q]=\sum_{\overline{g}\in P\backslash G/Q}\bigl[G/P\cap {}^gQ\bigr]
\]
\end{lemma}

\begin{proof}
The following is Nicolas Lemonine's proof in \cite{A9}. Each couple $(g_1P,g_2Q)\in (G/P)\times(G/Q)$ is in the same orbit as $(P,g_1^{-1}g_2Q)$ for the diagonal action of $G$. Moreover, $(P,gQ)$ and $(P,g'Q)$ are in the same orbit if and only if there exists $\widetilde{g}\in G$ such that $\widetilde{g}P=P$ and $g'^{-1}\widetilde{g}g\in Q$, i.e.\ if and only if $g\in Pg'Q$. To conclude, when $g\in G$ is fixed we can compute the stabilizer subgroup of $(P,gQ)$:
\begin{align*}
\operatorname{Stab}_G\bigl((P,gQ)\bigr)
&=\{g'\in G\mid g'\cdot (P,gQ)=(P,gQ)\}\\
&=\{g'\in G\mid g'\in P \text{ and } g^{-1}g'g\in Q\}\\
&=P\cap {}^gQ
\end{align*}
\end{proof}

We can check that:

\begin{lemma}
Let\ $G$ be a finite group, then the elements of\ $B(G)$ satisfy:
\par (1)\ $[G/P_{i}][G/P_{j}]=[G/P_{i} \times G/P_{j}]=[G/P_{j} \times G/P_{i}]=[G/P_{j}][G/P_{i}]$, where\ $P_{i},P_{j}\leq G$;
\par (2)\ $[G/P_{i} \times G/1]=|G/P_{i}|[G/1]$,\ $P_{i}\leq G$;
\par (3)\ $[G/P_{i} \times G/G]=[G/G \times G/P_{i}]=[G/P_{i}]$,\ $P_{i}\leq G$.
\end{lemma}

Now we recall the description of the Burnside ring $B(G)$ in terms of fixed points of $G$-sets.

\begin{definition}\cite[Definition 1.2]{A9}\label{definition 2.4}
For each subgroup $P \leq G$, we define a $\mathbb{Z}$-linear form $\Phi_P \colon B(G) \to \mathbb{Z}$ by setting
\[
\Phi_P([X]) := |X^P|
\]
for each isomorphism class of $G$-sets $[X]$, where $X$ is any representative of the class, and where $X^P$ denotes the set of fixed points for the $P$-action on $X$.
\end{definition}

The value of $\Phi_P$ on transitive $G$-sets is given by the following classical formula (see Section 3 in~\cite{A2}).

\begin{lemma}\label{lemma 2.5}\cite[Section 3]{A2}\cite[Lemma 1.3]{A9}
Let $P$ and $Q$ be subgroups of $G$. Then we have
\[
\Phi_P\bigl([G/Q]\bigr)=\frac{\bigl|N_G(P,Q)\bigr|}{|Q|}
\]
where $N_G(P,Q) := \{g \in G \mid {}^gP \leq Q\}$.
\end{lemma}

\begin{proof}
Nicolas Lemonine gave a short proof in \cite[Lemma 1.3]{A9}.
\end{proof}

\begin{lemma}
If\ $[X]$ is a unit of\ $B(G)$, then\ $[X]^{2}=1$.
\end{lemma}

\begin{proof}
Lemma $2.1$ can be interpreted as follows: each subgroup $H$ of $G$ defines a ring homomorphism
\[
\phi_H^G \colon B(G) \to \mathbb{Z}
\]
by $\phi_H^G(X) = |X^H|$. The kernel of $\phi_H^G$ is a prime ideal, since $\mathbb{Z}$ is an integral domain, and the intersection of all those kernels for subgroups $H$ of $G$ is zero.

Since $\phi_H^G = \phi_K^G$ if $H$ and $K$ are conjugate in $G$, it follows that the product map
\[
\Phi = \prod_{H\in Cl(G)} \phi_H^G \colon B(G) \to \prod_{H\in Cl(G)} \mathbb{Z},\ [X]\mapsto (|X^{H}|)_{H\in Cl(G)}
\]
is injective.

Since \(B(G)^\times \hookrightarrow \left( \prod_{H \leq G} \mathbb{Z} \right)^\times\), and
\[
\left( \prod_{H \leq G} \mathbb{Z} \right)^\times =  \prod_{H \leq G} \{\pm 1\},
\]
it follows that \(B(G)^\times\) is a finite elementary abelian \(2\)-group. This completes the proof.
\end{proof}

\subsection{The Burnside ring of a fusion system }
Now we will consider a finite $p$-group $S$ and fusion system $\mathcal{F}$ over $S$. We will denote the set of morphisms from
$P$ to $Q$ in $\mathcal{F}$ by $\mathcal{F}(P, Q)$. The preorder $\preceq_{\mathcal{F}}$ on subgroups of $S$ denotes $\mathcal{F}$-subconjugation, meaning that
\[
P \lesssim_{\mathcal{F}} Q \iff \exists \varphi \in \mathcal{F}(P,Q) \iff \exists \varphi \in \mathcal{F}(P,S),\ \varphi(P) \leq Q.
\]
When $P \lesssim_{\mathcal{F}} Q$ and $Q \lesssim_{\mathcal{F}} P$, the subgroups $P$ and $Q$ are $\mathcal{F}$-isomorphic (or $\mathcal{F}$-conjugate), and we denote this by $P \sim_{\mathcal{F}} Q$. We denote the set of $\mathcal{F}$-isomorphism classes of subgroups of $S$ by $Cl(\mathcal{F})$. The class of a subgroup $P\leq S$ in $Cl(\mathcal{F})$ will be denoted by $[P]_{\mathcal{F}}$.

\begin{definition}
An element $[X] \in B(S)$ is said to be $\mathcal{F}$-stable if for every subgroup $P \leq S$ and every morphism $\varphi \in \mathcal{F}(P,S)$, we have
\[
r_\varphi([X]) = r_{i_{P,S}}([X]) \quad \text{in } B(P),
\]
where $i_{P,S}$ is the inclusion $P \hookrightarrow S$.
\end{definition}

The following characterisation will be helpful. A version of this lemma was included in the Ph.D.\ thesis of Gelvin \cite[Proposition 3.2.3]{A7}, and previously the lemma has at least been implicitly used by Broto, Levi and Oliver \cite[Proof of Proposition 5.5]{A6}. A special case of the lemma, for bisets, was also proved by Ragnarsson and Stancu \cite[Lemma 4.8, parts (b) and (c)]{A8}. Nicolas applied this lemma in~\cite{A9}.  For the proof, we refer the reader to \cite{A2}.

\begin{lemma}\cite[Proposition 1.7]{A9}\cite[Lemma 4.1]{A2} \label{prop:Prop. 4.1}
Let $[X]$ be an element of $B(S)$. The following are equivalent:
\begin{enumerate}
    \item[(1)] $[X]$ is $\mathcal{F}$-stable.
    \item[(2)] $\Phi_P([X]) = \Phi_Q([X])$ whenever $P$ and $Q$ are $\mathcal{F}$-isomorphic subgroups of $S$.
    \end{enumerate}
\end{lemma}

To define the Burnside ring of a fusion system $\mathcal{F}$, we have two possibilities. We can consider the semiring of all the $\mathcal{F}$-stable $S$-sets and take the Grothendieck group of this. Alternatively, we can first take the Grothendieck group for all $S$-sets to get the Burnside ring of $S$, and then afterwards we consider the subring herein consisting of all the $\mathcal{F}$-stable elements. The following proposition implies that the two definitions coincide for saturated fusion systems.

\begin{lemma}\cite[Proposition 4.4]{A2} \label{prop:Prop. 4.4}
Let $\mathcal{F}$ be a fusion system over a $p$-group $S$, and consider the subsemiring $B_{+}(\mathcal{F})$ of $\mathcal{F}$-stable $S$-sets in the semiring $B_{+}(S)$ of finite $S$-sets. This inclusion induces a ring homomorphism from the Grothendieck group of $B_{+}(\mathcal{F})$ to the Burnside ring $B(S)$, which is injective.

If $\mathcal{F}$ is saturated, then the image of the homomorphism is the subring of $B(S)$ consisting of the $\mathcal{F}$-stable elements.
\end{lemma}

\begin{definition}\cite[Definition 1.8]{A9}\cite[Definition 4.5]{A2} \label{def:Def. 4.5}
The set of $\mathcal{F}$-stable elements is a subring of $B(S)$, which we call the Burnside ring of $\mathcal{F}$ and denote by $B(\mathcal{F})$, i.e.
\[
B(\mathcal F)
=
\left\{
[X] \in B(S)
\;\middle|\;
\Phi_P([X]) = \Phi_{\varphi(P)}([X])
\text{ for all }P\le S
\text{ and all }\varphi\in\operatorname{Hom}_{\mathcal F}(P,S)
\right\}.
\]
\end{definition}

\begin{lemma}
If\ $[X]$ is a unit in\ $B(\mathcal{F})$, then\ $[X]^{2}=1$.
\end{lemma}
\begin{proof}
For any $p$-group $S$, the unit group $B(S)^{\times}$ is an elementary abelian $2$-group. By Lemma \ref{prop:Prop. 4.4} and Definition \ref{def:Def. 4.5}, $B(\mathcal{F})^{\times}\leq B(S)^{\times}$ is also an elementary abelian $2$-group, which shows that $[X]^{2}=1$.
\end{proof}

Here we describe particular units of $B(S)^{\times}$ indexed by the maximal subgroups of $S$. These units are constructed in such a way that it is easy tell if they are contained in $B(\mathcal{F})$ for any $\mathcal{F}$ on $S$.

\begin{lemma}\cite[Proposition 8.1]{A1} \label{prop:Prop. 8.1}
Let $S$ be a $2$-group. For every maximal subgroup $M < S$, there exists a unit $1-[S/M]\in B(S)^{\times}$ where $1:= [S/S]$, such that, for $P \leq S$, one has
$$\Phi_P(1-[S/M]) = \begin{cases} 1, & PM = S, \\ -1, & PM = M. \end{cases}$$
\end{lemma}

\begin{proof}
This Lemma has already been proved by Barsotti and Carman in \cite[Proposition 8.1]{A1}; a short proof is given below.

Since $S$ is a $2$-group and $M$ is a maximal subgroup of $S$, $|S:M|=2$ and we can write $S/M = \{M, g_{0}M \}$ for some $g_{0}\in S\backslash M$. Then
\begin{equation}\notag
\begin{aligned}
S/M \times S/M=\{(M,M), (M,g_{0}M),(g_{0}M,M),(g_{0}M,g_{0}M) \} = (M,M)^{S}\sqcup (M,g_{0}M)^{S}.
\end{aligned}
\end{equation}
Note that
$$(M,M)^{S}=\{(M,M),(g_{0}M,g_{0}M) \},$$
$$(M,g_{0}M)^{S}=\{(M,g_{0}M),(M,g_{0}M) \},$$
$$Stab_{S}((M, M))=\{g\in S | g(M, M)=(M, M) \} = M ,$$
$$Stab_{S}((M, g_{0}M))=\{g\in S | g(M, g_{0}M)=(M, g_{0}M) \} = M ,$$
thus $(M,M)^{S}\cong S/Stab_{S}((M,M))=S/M$ and $(M, g_{0}M)^{S}\cong S/Stab_{S}((M, g_{0}M))=S/M$, which shows that $[S/M] [S/M]= [S/M \times S/M]=2[S/M] $. Hence $(1-[S/M])^{2} = 1 - 2[S/M] + [S/M]^{2} = 1$, i.e. $1-[S/M]\in B(S)^{\times}$.

For any $P \leq S$, since $M \lhd S $, the action of $P$ on $S/M$ is trivial if and only if $P \leq M$, which shows that $\Phi_P(1-[S/M]) = \begin{cases} 1, & PM = S, \\ -1, & PM = M. \end{cases}$ by Lemma \ref{lemma 2.5}.
\end{proof}

\subsection{The $1$-th cohomology group of the fusion system}\label{sec:The cohomology of the fusion system}
The following introduces some basic background and important lemmas needed in our proof, with reference to \cite{A6,Sh,W}.

\begin{definition}\cite[Definition 1.2.4]{Sh}\cite[Definition 6.4.1]{W} \label{def:Definition 6.4.1}
A $1$-cocycle (or crossed homomorphism) of group $G$ in a left $G$-module $A$ is a set map $f: G \rightarrow A$ satisfying
$$f(gh) = f(g) + g\cdot f(h)$$ for all $g,h \in G$. The family $Z^{1}(G,A)$ of all $1$-cocycle is an abelian group in an obvious way $(f_{1} + f_{2})(g)= f_{1}(g) + f_{2}(g)$.
\end{definition}
Note that if $A$ is a trivial $G$-module, $Z^{1}(G,A)= Hom(G,A)$.

\begin{definition}\cite[Definition 1.2.4]{Sh}\cite[Definition 6.4.2]{W} \label{def:Definition 6.4.2}
If $a\in A$, define $f_{a}(g) = g\cdot a - a$; $f_{a}$ is a $1$-cocycle because
$$f_{a}(gh)=(gh)\cdot a -g\cdot a + g\cdot a- a=g\cdot f_{a}(h) + f_{a}(g).$$
The $f_{a}$ are called the $1$-coboundaries of $G$ in $A$. Since $f_{a}+f_{b}=f_{a+b}$ for all $a,b\in A$, the set $B^{1}(G,A)$ of $1$-coboundaries forms a subgroup of $Z^{1}(G,A)$.
\end{definition}

\begin{definition}\cite[Definition 1.2.5, Lemma 1.2.6]{Sh}\cite[Theorem 6.4.5]{W} \label{def:Definition 1.2.5}
We define the $1$-th cohomology group of $G$ with coefficients in $A$ to be $H^{1}(G,A)=Z^{1}(G,A)/B^{1}(G,A)$.
\end{definition}

\begin{lemma}\cite[Lemma 1.2.6]{Sh}\cite[Corollary 6.4.6]{W} \label{Lemma 1.2.6}
If $A$ is a trivial $G$-module, then $H^{1}(G,A)=Z^{1}(G,A)/B^{1}(G,A)\cong Z^{1}(G,A) \cong \operatorname{Hom}(G,A)$.
\end{lemma}

\begin{definition}
Let $p$ be a prime, $S$ be a finite $p$-group, and regard $\mathbb{F}_{p}$ as a trivial $S$-module. The $1$-th cohomology group $H^1(S;\mathbb F_p)$ of $S$ with coefficients in $\mathbb{F}_{p}$  is defined as
\[
H^1(S;\mathbb F_p)=Z^1(S;\mathbb F_p)/B^1(S;\mathbb F_p).
\]
\end{definition}
Since $S$ acts trivially on $\mathbb{F}_{p}$, the $1$-cocycle condition becomes $f(gh)=f(g)+f(h)$ for all $g,h\in S$, while every
\(1\)-coboundary is zero. Hence there is a natural isomorphism
\[
H^1(S;\mathbb F_p)\cong \operatorname{Hom}(S,\mathbb F_p).
\]

\begin{lemma}\label{lemma:lemma Phi}
Let $p$ be a prime, let $S$ be a finite $p$-group, and regard $\mathbb{F}_{p}$ as a trivial $S$-module. Then $H^1(S;\mathbb F_p)\cong \operatorname{Hom}(S,\mathbb F_p) \cong \operatorname{Hom}\bigl(S/\Phi(S),\mathbb F_p\bigr) \cong \operatorname{Hom}_{\mathbb F_p}\bigl(S/\Phi(S),\mathbb F_p\bigr)\cong S/\Phi(S)$, where $\Phi(S)$ denotes the Frattini subgroup of $S$.
\end{lemma}
\begin{proof}
Since $\mathbb{F}_{p}^{+}$ is an abelian group, $[S,S]\leq Ker(f)$ for every $f\in \operatorname{Hom}(S,\mathbb{F}_{p})$. It is clear that $f(g^{p})=pf(g)=0$ for all $f\in \operatorname{Hom}(S,\mathbb{F}_{p})$ and $g\in S$, then $S^{p}=\langle g^{p}| g\in S\rangle \leq Ker(f)$ for every $f\in \operatorname{Hom}(S,\mathbb{F}_{p})$. Thus $\Phi(S)= S^{p}[S,S]\leq Ker(f)$. For every $f\in \operatorname{Hom}(S,\mathbb{F}_{p})$, we can write $f=\overline{f}\circ \pi$ for unique $\overline{f}\in \operatorname{Hom}\bigl(S/\Phi(S),\mathbb F_p\bigr)$, where $\pi:S\rightarrow S/\Phi(S)$.

Then $H^1(S;\mathbb F_p)\cong \operatorname{Hom}(S,\mathbb F_p) \cong \operatorname{Hom}\bigl(S/\Phi(S),\mathbb F_p\bigr)$.

Since $V=S/\Phi(S)$ is an elementary abelian p-group, it can naturally be regarded as a finite-dimensional vector space over $\mathbb F_p$. Let $\dim_{\mathbb F_p}V=n$, then $V\cong (\mathbb F_p)^n$. Its dual space $V^{*}=\operatorname{Hom}_{\mathbb F_p}(V,\mathbb F_p)$ is exactly $\operatorname{Hom}(V,\mathbb F_p)$, and $V^{*}\cong V$.

Thus $H^1(S;\mathbb F_p)\cong \operatorname{Hom}(S,\mathbb F_p) \cong \operatorname{Hom}\bigl(S/\Phi(S),\mathbb F_p\bigr) \cong \operatorname{Hom}_{\mathbb F_p}\bigl(S/\Phi(S),\mathbb F_p\bigr)\cong S/\Phi(S)$.
\end{proof}

\begin{remark}\label{remark:Remark S-p=2}
In particular, for $p=2$,
\[
H^1(S;\mathbb F_2)\cong \operatorname{Hom}(S,\mathbb F_2)\cong \operatorname{Hom}\bigl(S/\Phi(S),\mathbb F_2\bigr).
\]
\end{remark}

\begin{definition}\cite[Proposition 5.5, Theorem 5.8]{A6}\cite[\S2]{DGPS}\cite[Definition 2.3]{Mo} \label{def:Definition 2.3}
Let $\mathcal F$ be a saturated fusion system on a finite $p$-group $S$. The $1$-th cohomology of $\mathcal{F}$ is defined as the subring of $\mathcal F$-stable elements of $H^1(S;\mathbb F_p)$:
\[
H^1(\mathcal F;\mathbb F_p)
=
\left\{
\chi \in H^1(S;\mathbb F_p)
\;\middle|\;
\varphi^*(\chi)=\operatorname{res}_P^S(\chi)
\text{ for all }P\le S
\text{ and all }\varphi\in\operatorname{Hom}_{\mathcal F}(P,S)
\right\},
\]
where $\operatorname{res}_P^S$ is the restriction map and $\varphi^*$ is restriction along $\varphi$. By $H^1(S;\mathbb F_p)=\operatorname{Hom}(S,\mathbb F_p)$, the above condition is equivalent to
\[
\chi(\varphi(x))=\chi(x)
\qquad
\text{for all }P\le S,\ \varphi\in\operatorname{Hom}_{\mathcal F}(P,S),\ x\in P.
\]
Then we can write
\[
H^1(\mathcal F;\mathbb F_p)
=
\left\{
\chi \in \operatorname{Hom}(S;\mathbb F_p)
\;\middle|\;
\chi|_{P} = \chi \circ \varphi
\text{ for all }P\le S
\text{ and all }\varphi\in\operatorname{Hom}_{\mathcal F}(P,S)
\right\}.
\]
\end{definition}
The following lemma is a proof tool needed in this paper; for its detailed statement and proof, see \cite{DGPS}. We will give a brief proof later.

\begin{lemma}[{\cite[Corollary~1.2]{DGPS}}] \label{thm:DGPS-cor.1.2}
Let $\mathcal F$ be a saturated fusion system on a finite $p$-group $S$. If the restriction map
\[
H^1(\mathcal F;\mathbb{F}_p)\longrightarrow H^1(S;\mathbb F_p)
\]
is an isomorphism, then $\mathcal F=\mathcal F_S(S)$, the trivial fusion system on $S$.
\end{lemma}

To prove Lemma \ref{thm:DGPS-cor.1.2}, we need to introduce some definitions and theorems.

\begin{definition}\cite[Definition 1.1]{DGPS} \label{def:Def. 1.1}
Let $\mathcal{F}$ be a saturated fusion system on a finite $p$-group $S$.\\
(1) $O_{\mathcal{F}}^p(S) = \langle [P, O^p(\mathrm{Aut}_{\mathcal{F}}(P))] \mid P \leqslant S \rangle$ (the hyperfocal subgroup of $\mathcal{F}$), where $O^p(\mathrm{Aut}_{\mathcal{F}}(P))= \bigcap \{N\unlhd \mathrm{Aut}_{\mathcal{F}}(P)| \mathrm{Aut}_{\mathcal{F}}(P)/N \text{ are}\  p-\text{group}\} $;\\
(2) $A_{\mathcal{F}}^p(S) = [S, \mathcal{F}] = \langle [P, \mathrm{Aut}_{\mathcal{F}}(P)] \mid P \leqslant S \rangle = \langle x^{-1}\alpha(x) \mid x\in P, P \leqslant S, \alpha \in \mathrm{Aut}_{\mathcal{F}}(P) \rangle$ (the focal subgroup of $\mathcal{F}$); and\\
(3) $E_{\mathcal{F}}^p(S) = \Phi(S)[S, \mathcal{F}] = \Phi(S)O_{\mathcal{F}}^p(S)$ (the elementary focal subgroup of $\mathcal{F}$).
\end{definition}

\begin{lemma}\cite[Theorem T]{DGPS} \label{thm:DGPS-Theorem T}
(Tate's theorem for fusion systems). Let $\mathcal{F}$ be a saturated fusion system on a finite $p$-group $S$ and let $\mathcal{H}$ be a saturated subsystem of $\mathcal{F}$ on $S$. The following are equivalent.\\
(1) $E_{\mathcal{F}}^p(S) = E_{\mathcal{H}}^p(S)$;\\
(2) $A_{\mathcal{F}}^p(S) = A_{\mathcal{H}}^p(S)$;\\
(3) $O_{\mathcal{F}}^p(S) = O_{\mathcal{H}}^p(S)$.
\end{lemma}
\begin{proof}
This is one of the main results of the paper \cite{DGPS} by A. Díaz, A. Glesser, S. Park and R. Stancu, for its detailed proof, see \cite{DGPS}.
\end{proof}

\begin{proof}[Proof of Lemma \ref{thm:DGPS-cor.1.2}]
We shall prove first that $H^1(\mathcal{F}; \mathbb{F}_p) \cong \text{Hom}(S/[S, \mathcal{F}], \mathbb{F}_p)$. By Definition \ref{def:Def. 1.1} $(2)$, for any $x^{-1}\alpha(x)\in A_{\mathcal{F}}^p(S) = [S, \mathcal{F}] = \langle x^{-1}\alpha(x) \mid x\in P, P \leqslant S, \alpha \in \mathrm{Aut}_{\mathcal{F}}(P) \rangle$ and $\chi \in H^1(\mathcal{F}; \mathbb{F}_p),$ we have $\chi(x)= \chi(\alpha(x))$ for all $x\in P$, which shows that $\chi(x^{-1}\alpha(x))= 0 $ for all $x^{-1}\alpha(x)\in A_{\mathcal{F}}^p(S)= [S, \mathcal{F}]$. Then we can check that $H^1(\mathcal{F}; \mathbb{F}_p) \rightarrow \text{Hom}(S/[S, \mathcal{F}], \mathbb{F}_p), \chi \mapsto \overline{\chi}$ is an isomorphism, where $\chi = \overline{\chi}\circ \pi$ and $\pi:S\rightarrow S/[S,\mathcal{F}]$.

By Lemma \ref{lemma:lemma Phi}, $\text{Hom}(S/[S,\mathcal{F}], \mathbb{F}_p) \cong \text{Hom}(S/\Phi(S)[S,\mathcal{F}],\mathbb{F}_p) \cong S/\Phi(S)[S,\mathcal{F}] = S/E_{\mathcal{F}}^p(S).$ Indeed, $\Phi(S/[S,\mathcal{F}])= \Phi(S)[S,\mathcal{F}]/[S,\mathcal{F}]$, then $S/[S,\mathcal{F}] \Big/ \Phi(S/[S,\mathcal{F}]) = S/[S,\mathcal{F}] \Big/ \Phi(S)[S,\mathcal{F}]/[S,\mathcal{F}]\cong S/\Phi(S)[S,\mathcal{F}]$, which shows that $\text{Hom}(S/[S,\mathcal{F}], \mathbb{F}_p) \cong \text{Hom}(S/\Phi(S)[S,\mathcal{F}],\mathbb{F}_p)$.

Clearly, the trivial fusion system $\mathcal{H} = \mathcal{F}_{S}(S)$ is a saturated subsystem of $\mathcal{F}$ on $S$, we obtain $$H^1(\mathcal{F}_S(S); \mathbb{F}_p) = \text{Hom}(S/[S, \mathcal{F}_S(S)], \mathbb{F}_p)= \text{Hom}(S/[S,S], \mathbb{F}_p)\cong \text{Hom}(S, \mathbb{F}_p) \cong H^1(S; \mathbb{F}_p),$$
where $[S, \mathcal{F}_S(S)] = [S,S]$. Similar to Lemma \ref{lemma:lemma Phi}, $\text{Hom}(S/[S,S], \mathbb{F}_p) \cong S/\Phi(S) = S/\Phi(S)[S,S]= S\big/E_{\mathcal{F}_{S}(S)}^p(S)$.

From the hypothesis, $H^1(\mathcal F;\mathbb F_p)\cong H^1(S;\mathbb F_p)$, we get
\begin{align*}
S/E_{\mathcal{F}_{S}(S)}^{p}(S)= S/\Phi(S) \cong H^1(\mathcal{F}_S(S); \mathbb{F}_p) \cong H^1(S; \mathbb{F}_p) \cong H^1(\mathcal{F}; \mathbb{F}_p)\cong S/\Phi(S)[S,\mathcal{F}] = S/E_{\mathcal{F}}^p(S),
\end{align*}
then $E_{\mathcal{F}}^p(S)=E_{\mathcal{F}_{S}(S)}^{p}(S)$. By Lemma \ref{thm:DGPS-Theorem T}, $O_{\mathcal{F}}^p(S) = O_{\mathcal{F}_{S}(S)}^p(S)$. Since $\mathrm{Aut}_{\mathcal{F}_{S}(S)}(P) = \mathrm{Aut}_{S}(P)$ are $p$-groups for all $P\leq S$, $O^{p}(\mathrm{Aut}_{S}(P))=1$, which shows that $O_{\mathcal{F}}^p(S) = O_{\mathcal{F}_{S}(S)}^p(S) = 1 $. For any $P\leq S$ and $\beta\in O^p(\mathrm{Aut}_{\mathcal{F}}(P))$, we have $x^{-1}\beta(x)\in O_{\mathcal{F}}^p(S)=1$ for all $x\in P$, then $\beta(x)=x
$ for all $x\in P$ and $\beta$ is the identity map on $P$. Hence $O^p(\mathrm{Aut}_{\mathcal{F}}(P))=1$ for all $P\leq S$, which shows that $\mathrm{Aut}_{\mathcal{F}}(P)$ is a $p$-group for all $P\leq S$.

Since $\mathcal{F}$ is a saturated fusion system on a finite $p$-group $S$, $\mathrm{Aut}_{\mathcal{F}}(P) = \mathrm{Aut}_{S}(P)$ and $\mathcal{F} = \mathcal{F}_{S}(S)$.
\end{proof}

\begin{remark}\label{remark:Remark F-p=2}
For \(p=2\),
\[
H^1(\mathcal F;\mathbb F_2)
=
\left\{
\chi \in \operatorname{Hom}(S;\mathbb F_2)
\;\middle|\;
\chi|_{P} = \chi \circ \varphi
\text{ for all }P\le S
\text{ and all }\varphi \in \operatorname{Hom}_{\mathcal F}(P,S)
\right\}.
\]
\end{remark}
In particular, for \(p=2\), if
\[
H^1(\mathcal F;\mathbb F_2)\longrightarrow H^1(S;\mathbb F_2)
\]
is an isomorphism, then \(\mathcal F=\mathcal F_S(S)\).

\section{Proof of Theorem 1.2}
We now prove the main result of this paper. If $S=1$, then there is the only trivial fusion system on \(S\) and $B(1)\cong \mathbb{Z}, B(1)^{\times} \cong \{ 1, -1 \}$.

\begin{proof}[Proof of Lemma \ref{theorem:A}]
$(1) \Rightarrow (2)$ is trivial. By Lemma \ref{prop:Prop. 8.1}, we have $(2) \Rightarrow (3)$.

$(3) \Rightarrow (1)$ For any $0 \neq \tau \in \operatorname{Hom}(S,\mathbb F_2) \cong H^1(S;\mathbb F_2)$(Remark \ref{remark:Remark S-p=2}
), $P \leq S$ and $x \in P$, $\tau$ is surjective and $M = Ker(\tau)$ is a maximal subgroup of $S$ since \(S\) is a finite \(2\)-group. From condition $(3)$, $1-[S/M] \in B(\mathcal F)^\times$, i.e. $1-[S/M]$ is $\mathcal F$-stable. By Lemma \ref{prop:Prop. 4.1} and Definition \ref{def:Def. 4.5}, $\Phi_P(1-[S/M]) = \Phi_{\varphi(P)}(1-[S/M])$ for any $\varphi \in \operatorname{Hom}_{\mathcal F}(P,S)$, where the  $\Phi_P$ and $\Phi_{\varphi(P)}$ are defined in Definition \ref{definition 2.4}.

Let $\langle x \rangle$ be the cyclic group generated by $x$, we have $\varphi|_{\langle x\rangle} : \langle x\rangle \rightarrow \langle \varphi(x)\rangle$, where $\varphi|_{\langle x\rangle} \in \operatorname{Hom}_{\mathcal F}(\langle x\rangle,S)$. Since $1-[S/M] \in B(\mathcal F)^\times$, we also have $\Phi_{\langle x\rangle}(1-[S/M]) = \Phi_{\langle \varphi(x)\rangle}(1-[S/M])$.
By Lemma \ref{prop:Prop. 8.1},
$$\Phi_{\langle x\rangle}(1-[S/M]) = \begin{cases} 1, & \langle x\rangle M = S, \\ -1, & \langle x\rangle M = M. \end{cases}$$
and
$$\Phi_{\langle \varphi(x)\rangle}(1-[S/M]) = \begin{cases} 1, & \langle \varphi(x)\rangle M = S, \\ -1, & \langle \varphi(x)\rangle M = M. \end{cases}$$
Then $x \in M $ if and only if $\varphi(x) \in M$.

Since $\tau(g) = \begin{cases} 0, & g\in M \\ 1, & g \notin M \end{cases}$ for any $g \in S$, $\tau(\varphi(x)) = \begin{cases} 0, & \varphi(x)\in M \\ 1, & \varphi(x) \notin M \end{cases}$, which shows that $\tau(\varphi(x)) = \tau(x) = \tau|_{P}(x)$ for any $x\in P$. Hence $\tau|_{P} = \tau \circ \varphi $ for any $P \leq S$ and $\varphi \in \operatorname{Hom}_{\mathcal F}(P,S)$.

By Remark \ref{remark:Remark F-p=2},
\[
H^1(\mathcal F;\mathbb F_2)
=
\left\{
\chi \in \operatorname{Hom}(S;\mathbb F_2)
\;\middle|\;
\chi|_{P} = \chi \circ \varphi
\text{ for all }P\le S
\text{ and all }\varphi \in \operatorname{Hom}_{\mathcal F}(P,S)
\right\},
\]
then $\tau \in H^1(\mathcal F;\mathbb F_2)$, which shows that $H^1(S;\mathbb F_2)  \cong \operatorname{Hom}(S,\mathbb F_2) \subseteq H^1(\mathcal F;\mathbb F_2)$. Clearly, $ H^1(\mathcal F;\mathbb F_2) \subseteq H^1(S;\mathbb F_2) $, hence $ H^1(\mathcal F;\mathbb F_2) = H^1(S;\mathbb F_2)$. Thus the restriction map $H^1(\mathcal F;\mathbb F_2) \rightarrow H^1(S;\mathbb F_2)$ is an isomorphism, which shows that $ \mathcal F = \mathcal F_S(S)$ by Lemma \ref{thm:DGPS-cor.1.2}.
\end{proof}

\end{document}